\documentclass[12pt,psamsfonts]{amsart}

\newtheorem{theorem}{Theorem}[section]

\usepackage[all]{xy}                        

\usepackage[german,english]{babel}          
\usepackage{amssymb}                        
\usepackage{latexsym}
\usepackage{xspace}                         
\usepackage{stmaryrd}                       
\usepackage{calc}
\usepackage{ifthen}
\usepackage[margin=1.3in]{geometry}

\numberwithin{equation}{section}

\theoremstyle{plain}

\newtheorem*{theorem*}{Theorem}
\newtheorem{proposition}[theorem]{Proposition}
\newtheorem{corollary}[theorem]{Corollary}
\newtheorem{lemma}[theorem]{Lemma}

\theoremstyle{definition}
\newtheorem{definition}[theorem]{Definition}
\newtheorem{example}[theorem]{Example}

\newtheorem{question}[theorem]{Question}

\newtheorem{remark}[theorem]{Remark}

\theoremstyle{remark}

\newboolean{book} 
\setboolean{book}{false}

\newcommand{\new}{\newcommand}

\newcommand{\ren}{\renewcommand}

\renewcommand{\AA}{\mathbb{A}}

\newcommand{\CC}{\mathbb{C}}

\newcommand{\NN}{\mathbb{N}}

\newcommand{\PP}{\mathbb{P}}

\newcommand{\RR}{\mathbb{R}}

\newcommand{\ZZ}{\mathbb{Z}}

\newcommand  {\Char}    {\operatorname{char}}

\renewcommand{\O}       {\mathcal{O}}

\newcommand  {\rad}      {\operatorname{rad}}

\newcommand  {\reg}     {\operatorname{reg}}

\newcommand  {\Spec}    {\operatorname{Spec}}

\newcommand{\comdots}{ , \ldots , }
\newcommand{\komdots}{ , \ldots , }
\newcommand{\plusdots}{ + \ldots + }

\newcommand{\timesdots}{ \times \ldots \times }

\renewcommand{\bar}[1]{\overline{#1}}

\newcommand{\longto}{\longrightarrow}

\renewcommand{\phi}{\varphi}        
\renewcommand{\epsilon}{\varepsilon}

\newcommand{\tensor}{\otimes}         

\newcommand{\set}[1]{\left\{#1\right\}}

\renewcommand{\to}[1][]{\xrightarrow{\ #1\ }}

\usepackage{amscd}
\usepackage{amssymb}

\newboolean{extra} 

\setboolean{extra}{false}

\ren{\sec}[1]{{ {#1}'}}

\newcommand{\elef}{{f}}
\newcommand{\eleg}{{g}}

\newcommand{\elel}{{l}}
\newcommand{\elem}{{m}}
\newcommand{\elen}{{n}}

\newcommand{\elet}{{t}}

\ren{\elef}{t} 

\ren{\eleg}{s}

\ren{\elel}{t} 

\ren{\elem}{s} 

\ren{\elen}{t}

\ren{\elet}{y}

\newcommand{\algforc}{B}
\newcommand{\forcalg}{\algforc}

\newcommand{\mata}{D}

\newcommand{\dvd}{V} 

\newcommand{\modul}{{M}}

\newcommand{\submod}{{N}}

\newcommand{\modra}{\mu}
\newcommand{\submodra}{\nu}

\newcommand{\spax}{X}

\newcommand{\spaxeqspecring}{\spax = \Spec \ring }

\newcommand{\ringr}{R}
\newcommand{\rings}{S}

\newcommand{\covmap}{\varphi}

\newcommand{\covmapcov}{\covmap_i: X_i \to X}

\newcommand{\synchronize}[2]{\renewcommand{#1}{{#2}}}

\newcommand{\ind}{\kappa}

\newcommand{\indi}{i}

\newcommand{\indj}{j}

\newcommand{\runi}{{i}}

\newcommand{\Indj}{{J}}

\newcommand{\realnum}{\alpha}

\newcommand{\realnumbeta}{\beta}
\newcommand{\realnumlambda}{\lambda}

\newcommand{\compnumt}{{t}}

\newcommand{\brackconvl}{{\{\!\{}} 

\newcommand{\brackconvr}{{\}\!\}}}

\newcommand{\compconj}{\bar} 

\newcommand{\ebd}{{k}} 

\newcommand{\spay}{Y}

\newcommand{\spav}{{V}}

\newcommand{\spac}{{C}}

\newcommand{\point}{{P}}

\newcommand{\numpoints}{{k}}

\newcommand{\sphere}{{S}}

\newcommand{\neigh}{{U}}

\newcommand{\cofu}{\psi}
\newcommand{\contfu}{q}
\newcommand{\contfuhi}{ \varphi}

\newcommand{\contmap}{{\varphi}}

\newcommand{\proba}{h} 

\newcommand{\sumnormquad}{{\normcomp \fuf_1 \normcomp ^2 \plusdots \normcomp \fuf_\numgen \normcomp^2}}

\new{\normcomp}{{|}}

\newcommand{\fubound}{{q}}

\newcommand{\group}{{G}}

\newcommand{\elgroup}{{\xi}}

\newcommand{\field}{{K}}

\newcommand{\fieldl}{{L}}

\newcommand{\fieldres}{{\kappa}}

\newcommand{\tupelpoint}{{\alpha}} 

\newcommand{\indpoint}{{\kappa}}

\newcommand{\unitroot}{{\zeta}} 

\newcommand{\ring}{{R}}

\newcommand{\ringsec}{{S}}

\newcommand{\ringt}{{T}}

\newcommand{\homtest}{{\varphi}}

\newcommand{\mortest}{{\psi}}

\newcommand{\homring}{{\varphi}}

\newcommand{\ringhom}{{\theta}}

\newcommand{\seminor}{{\rm sn}}

\newcommand{\fuf}{{f}}
\newcommand{\fug}{{g}}
\newcommand{\fuh}{{h}}

\newcommand{\fu}{{\fuf}}

\newcommand{\idgenfuf}{{\fuf_1 \comdots \fuf_\numgen}}

\newcommand{\runfuf}{{\fuf_1 \comdots \fuf_\numgen}}

\newcommand{\indfu}{{i}} 

\newcommand{\degfu}{{d}}

\newcommand{\runfu}{{ \fu_1 \comdots \fu_\numgen}}

\newcommand{\numvar}{m}

\newcommand{\var}{z}
\newcommand{\varx}{z}
\newcommand{\vary}{w}
\newcommand{\varz}{v}

\newcommand{\varxx}{x}
\newcommand{\varyy}{y}
\newcommand{\varzz}{z}

\newcommand{\indvar}{j}

\newcommand{\expox}{{r}}
\newcommand{\expoy}{{s}}

\newcommand{\expon}{{e}}

\newcommand{\powser}{G}

\newcommand{\dege}{d}

\newcommand{\degd}{d}
\newcommand{\degm}{{\delta}}

\newcommand{\dime}{n}
\newcommand{\mondeg}{d}
\newcommand{\expot}{\tau}

\newcommand{\expofu}{{n}}

\new{\expoideal}{{d}}

\new{\varcoef}{{t}}  

\new{\fucoef}{{g}}

\new{\indcoef}{{\nu}}

\new{\vargen}{{z}}

\newcommand{\schemetest}{{C}}

\newcommand{\schemesn}{{D}}

\newcommand{\numcomp}{{m}}

\newcommand{\ideal}{{I}}

\newcommand{\numgen}{{n}} 

\newcommand{\numfu}{{n}}

\newcommand{\indgen}{{i}}

\newcommand{\intclo}{\overline }

\newcommand{\indseq}{{n}}  

\new{\fuarb}{{g}}

\newcommand{\varfu}{q}

\newcommand{\vart}{{T}} 

\newcommand{\multexp}{{\gamma}}

\newcommand{\multexptau}{\tau}

\newcommand{\multexpsigma}{\gamma}

\newcommand{\setexp}{\Gamma}

\newcommand{\monom}{{\var^\multexp}}

\newcommand{\monomideal}{{\ideal}}

\newcommand{\coeflins}{{\delta}} 

\newcommand{\transp}{{t}} 

\newcommand{\valu}{{\nu}}

\newcommand{\axis}{{L}}

\newcommand{\indaxis}{{j}}

\newcommand{\indaxisrun}{{\iota}}

\newcommand{\numaxes}{{k}} 

\newcommand{\valaxis}{{\rm val}} 
\newcommand{\ordaxis  }{{\mu}} 

\newcommand{\varax}{{x}} 
\newcommand{\varaxv}{{y}}

\newcommand{\ello}{{\ell}}

\newcommand{\schemeaxes}{{C}}
\newcommand{\schemeaxesd}{{D}}

\newcommand{\indvaraxi}{i}

\newcommand{\indvaraxj}{j}

\newcommand{ \equataxes}{{\varax_\indvaraxi \varax_\indvaraxj, \indvaraxi \neq
\indvaraxj}}

\newcommand{\orig}{{P}}

\newcommand{\degaxis}{\delta}

\newcommand{\wsi}{{\rm wsi}}
\newcommand{\ax}{{\rm ax}}

\newcommand{\cont}{{\rm cont}}

\newcommand{\fupsi}{\psi}

\newcommand{\coefc}{c}

\newcommand{\polp}{p}

\ren{\set}{{S}}

\ren{\modra}{{m}}

\ren{\submodra}{{\numgen}}

\ren{\expofu}{{k}}

 \ren{\indaxisrun}{{i}}

\ren{\varfu}{{w}}

\ren{\indseq}{{m}}

\begin{document}

\title[Continuous solutions to algebraic forcing equations]
{Continuous solutions to algebraic forcing equations}

\author[Holger Brenner]{Holger Brenner}
\address{Department of Pure Mathematics, University of Sheffield,
Hicks Building, Houns\-field Road, Sheffield S3 7RH, United Kingdom}
\email{H.Brenner@sheffield.ac.uk}





\begin{abstract}
We ask for a given system of polynomials $\fuf_1 \comdots
\fuf_\numgen$ and $\fuf$ over the complex numbers $\CC$ when there
exist continuous functions $\contfu_1 \comdots \contfu_\numgen$ such
that $\contfu_1 \fuf_1 \plusdots \contfu_\numgen \fuf_\numgen =
\fuf$. This condition defines the continuous closure of an ideal. We
give inclusion criteria and exclusion results for this closure in
terms of the algebraically defined axes closure. Conjecturally,
continuous and axes closure are the same, and we prove this in the
monomial case.
\end{abstract}

\maketitle


\section*{Introduction}

Let $\ring=\CC[\var_1 \comdots \var_\numvar]$ be the polynomial ring
over the complex numbers. Every polynomial $\fuf \in \ring$ defines
a (holomorphic and) continuous function $\CC^\numvar \to \CC$. Let
$\ideal=(\fuf_1 \comdots \fuf_\numgen)$ be an ideal and $\fuf \in
\ideal$. The containment $\fuf \in \ideal$ means that there exist
polynomials $\fuarb_1 \comdots \fuarb_\numgen \in \ring$ such that
$\fuf=\fuarb_1 \fuf_1 \plusdots \fuarb_\numgen \fuf_\numgen$ holds.
If we allow $\fuarb_1 \comdots \fuarb_\numgen$ to be holomorphic
functions or formal power series, then we cannot express additional
polynomials as such a linear combination. On the other hand, if we
allow arbitrary functions $\CC^\numvar \to \CC$, then it is easy to
see that every polynomial inside the radical $\rad(\ideal)$ can be
expressed in this way.

In this paper we address the question under which condition we can
write  $\fuf = \contfu_1 \fuf_1 \plusdots \contfu_\numgen
\fuf_\numgen$ with continuous functions $\contfu_i :\CC^\numvar \to
\CC$. Putting it another  way: when does the so-called \emph{forcing
equation} $\fuf_1 \vart_1 \plusdots  \fuf_\numgen \vart_\numgen=
\fuf$ has a continuous solution? For a given ideal, the set of $\fuf
\in \ring$ such that there exists a continuous solution form an
ideal, which we denote by $\ideal ^\cont$ and which we call the
\emph{continuous closure} of $\ideal $ (Section \ref{contsection}).
For a principal ideal $\ideal =(\fug)$ we always have $\ideal
=\ideal^\cont$, since a rational function $\fuf/\fug$ has a
continuous extension to $\CC^\numvar $ only if $\fuf$ is a multiple
of $\fug$ in $\ring$. The easiest example for $\ideal \neq
\ideal^\cont$ is given by $ \varx^2 \vary^2 \not\in ( \varx^3,
\vary^3)$, but $ \varx^2 \vary^2 \in ( \varx^3, \vary^3)^\cont$ (we
will also see that $ \varx \vary \not\in ( \varx^2,
\vary^2)^\cont$).

Our main question is whether there exists an algebraic description
of the continuous closure. An easier task is to give good
algebraically defined approximations. For this we introduce the
\emph{axes closure} (Section \ref{axessection}), which is defined in
the following way: an element $\fuf$ in a commutative ring $\ring$
belongs to the ($\field$-)axes closure $\ideal ^\ax$ if and only if
for every ring homomorphism $\homtest: \ring \to \ringt$, where
$\Spec \ringt$ is a \emph{scheme of axes} over $\field$,
$\homtest(\fuf)$ belongs to the extended ideal $\homtest
(\ideal)\ringt$. A scheme of axes over $\field$ is given by a one
dimensional $\field$-algebra of finite type with normal components
meeting in one closed point $\point$ such that the completion (at
the meeting point) is isomorphic to $\fieldres(\point)[[\varax_1
\comdots \varax_\numaxes]]/(\equataxes)$ (Section
\ref{axesschemesection}).

The basic observation is that for such rings over $\field= \CC$ the
identity $\ideal =\ideal^\cont$ holds (Lemma \ref{axescont}),
therefore they serve as a category of testrings. It follows that
$\ideal ^\cont \subseteq \ideal^\ax$ (Corollary \ref{continaxes}),
which implies strong restrictions, since $\ideal ^\ax$ itself is
inside the weak subintegral closure (see
\cite{leahyvitulliweaksubintegralclosure},
\cite{reidvitullimonomial}) and in particular inside the integral
closure of the ideal. For example, $ \varx \vary \not\in ( \varx^2,
\vary^2)^\ax$, since restricted to the cross given by $( \varx+
\vary)( \varx- \vary)=0$ the element does not belong to the extended
ideal (on the two branches there exist solutions, but they do not
fit together), and hence $ \varx \vary \not\in ( \varx^2,
\vary^2)^\cont$ (but $\varx\vary \in \intclo{(\varx^2,\vary^2)}$,
its integral closure). We conjecture that $\ideal
^\cont=\ideal^\ax$, and we prove this in the case of a monomial
ideal.

To obtain inclusion results for the continuous closure one has to
construct continuous functions. Suppose that $\ideal =(\idgenfuf)$
is an ideal primary to the maximal ideal $(\var_1 \comdots
\var_\numvar)$ and $\fuf \in \ring$. In this case we can write
$\fuf= \sum_{\runi =1}^\numgen \frac{\fuf
\compconj{\fuf}_\indi}{\sumnormquad} \fuf_\indi$. Here $\frac{\fuf
\compconj{\fuf}_\indi}{\sumnormquad}$ is a continuous function on
$\CC^\numvar -\{0\} $, and the question is whether we can extend
this function continuously to the whole $\CC^\numvar$. This can be
done in the homogeneous situation under certain degree conditions
(Theorem \ref{degreecriterioncont}). We also show that $\fuf \in
\ideal^\cont$ under the condition that there exist numbers $\expofu
< \expoideal $ such that $\fuf^\expofu \in \ideal^\expoideal$
(Theorem \ref{powerscompare}, Corollary \ref{orderbiggeragain}).

For a monomial ideal $\monomideal=(\monom, \multexp \in \setexp)$,
$\setexp \subseteq \NN^\numvar$, it follows that the monomials
$\var^\expot$ with exponent $\expot$ in the interior of the convex
hull $\intclo{\setexp}$ of $\setexp$ in $\RR_+^\numvar$ belong to
the continuous closure (Theorem \ref{monomialaxescont}). The same
inclusion result holds for the axes closure over an arbitrary field
(Theorem \ref{monomialaxes}). For a monomial $\var^\expot$ with
exponent $\expot \in \intclo{\setexp} - {\intclo{\setexp}}^\circ$ on
the border we will show that $\var^\expot \not\in \ideal^\ax$
(unless $\var^\expot \in \ideal$) by reducing to the case where
$\ideal$ is generated by all monomials of degree $\degd$ with the
exception of one monomial $\var^\expot$ (Corollary
\ref{monomialbordernotaxes}). These results put together prove for a
monomial ideal that $\monomideal^\cont = \monomideal^\ax$ is the
monomial ideal $(\monom, \multexp \in \setexp \cup
{\intclo{\setexp}}^\circ))$ (Theorem \ref{monomialaxes} and Theorem
\ref{monomialaxescont}).

This work started in looking for a test category of rings for the
weak subintegral closure
(\cite{leahyvitulliweaksubintegralclosure},\cite{brennersemiintegraltest})
as presented by M. Vitulli at the conference on `Valuation Theory
and Integral Closures in Commutative Algebra', Ottawa, Canada, July
2006, so I would like to thank the organizers of this conference. I
thank T. Gaffney and M. Vitulli for discussions during the
conference and after. Furthermore I thank V. Bavula, D. Gepner, M.
Hochster, N. Jarvis, A. Kaid, M. Katzman, J. Manoharmayum, R.
Sanchez and R. Sharp for their interest and helpful remarks.

\section{The continuous closure}
\label{contsection}

\begin{definition}
Let $\ring$ be a finitely generated $\CC$-algebra, $\spax = \Spec
\ring $ and $\spax(\CC)$ the corresponding complex space with the
complex topology. Let $\ideal =(\fuf_1 \comdots \fuf_\numgen)$ be an
ideal in $\ring$ and $\fuf \in \ring$. We say that $\fuf \in
\ideal^\cont$, the \emph{continuous closure} of $\ideal $, if there
exist continuous functions $\contfu_\indgen: \spax(\CC) \to \CC$
such that $\fuf = \contfu_1 \fuf_1 \plusdots \contfu_\numgen
\fuf_\numgen$ as continuous functions on $\spax(\CC)$.
\end{definition}

Another way to think of the continuous closure is to look at the
ring homomorphism $\ring \to C_\CC (\spax(\CC))$ (which is an
inclusion if $\ring$ is reduced), where we denote by $C_\CC (\spax)$
the ring of continuous complex-valued functions on a topological
space $\spax$. Then $\ideal ^\cont$ is the contraction of the
extension of the ideal $\ideal $. In particular, $\ideal ^\cont$
does not depend on the choice of ideal generators.

\begin{remark}
It is also sometimes helpful to see the continuous closure in the
context of forcing algebras. This is the $\ring$-algebra $\forcalg =
\ring[\vart_1 \comdots \vart_\numgen]/( \fuf_1 \vart_1 \plusdots
\fuf_\numgen \vart_\numgen + \fuf)$. Set $\spay= \Spec \forcalg$.
Then $\fuf \in \ideal^\cont$ if and only if the natural projection
$\spay(\CC) \to \spax (\CC)$ has a continuous section $s: \spax(\CC)
\to \spay(\CC)$. This viewpoint allows us immediately to extend the
notion of continuous closure of an ideal to the continuous closure
of an $\ring$-submodule $\submod \subseteq \modul$ inside a finitely
generated $\ring$-module $\modul$. If $ \ring^\submodra
\stackrel{\mata}{\longrightarrow} \ring^\modra \longrightarrow
\modul/\submod \longrightarrow 0$ is a representation, where $\mata=
(\fuf_{\indi \indj})_{\indi \indj}$ is a $\modra \times \submodra
$-matrix with entries in $\ring$, and $\tilde{\elem}= (\elem_1
\comdots \elem_\modra) \in \ring^\modra$ denotes a representation of
$\elem \in \modul$, then $\elem \in \submod^\cont$ if and only if
the matrix equation $\mata (\contfu_1 \comdots \contfu_\submodra
)^\transp= (\elem_1 \comdots \elem_\modra)^\transp$ has a continuous
solution. However, we will not pursuit this generalization in this
paper. Neither will we be concerned with the Grothendieck topology
which arises by taking morphisms of finite type which allow a
continuous section as covers (see \cite{brennergrothendieck}).
\end{remark}

Our first inclusion result for the continuous closure is given by
the following theorem.

\begin{theorem}
\label{degreecriterioncont} Let $\ring$ be a finitely generated
standard-graded $\CC$-algebra and $\spax =\Spec \ring $,
$\spax(\CC)$ the corresponding complex analytic space. Let $\ideal
\subseteq \ring$ be a homogeneous $\ring_+$-primary ideal with
homogeneous ideal generators $\fuf_1 \comdots \fuf_\numgen$ of
degree $\degd_1 \comdots \degd_\numgen$. Let $\fuf$ be a homogeneous
element of degree $\degd > \max {\degd_\indi}$. Then there exist
continuous complex-valued functions $\contfu_\indi : \spax(\CC) \to
\CC$ such that $\fuf =  \contfu_1 \fuf_1 \plusdots \contfu_\numgen
\fuf_\numgen$ or in other words $\fuf \in \ideal^\cont$.
\end{theorem}
\begin{proof}
Since the only common zero of the $\fuf_\indi$ is the vertex (the
origin) $0 \in \spax$, they generate the unit ideal in
$C_\CC(\spax(\CC) -\{ 0 \} )$ (see the following Remark
\ref{contfuremark}). So there exist continuous functions
$\contfuhi_\indi :\spax(\CC) -\{ 0 \} \to \CC$ such that $\fuf
(\var) = \contfuhi_1(\var) \fuf_1(\var) \plusdots
\contfuhi_\numgen(\var) \fuf_\numgen(\var)$ for all $\var \in
\spax(\CC) - \{0\}$.

Consider a closed homogeneous embedding $\spax \subseteq
\AA_\CC^\ebd$ and $\spax (\CC) \subseteq \CC^\ebd$. Set $\sphere=
\spax (\CC) \cap \sphere^{2\ebd-1}$, where $\sphere^{2\ebd-1}=\{\var
\in \CC^\ebd : \normcomp \var\normcomp =1 \}$ is the real sphere.
Note that $\sphere$ is a compact closed subset of $\spax (\CC)$. We
restrict $\contfuhi_i$, $\indi=1 \comdots \numgen $, to $\sphere$ to
get continuous functions $\contfuhi_\indi :\sphere \to \CC$. The
idea is to extend these functions homogeneously to $\spax (\CC)
-\{0\}$ and to use the degree condition to show that it is even
extendable to the vertex.

Let $\proba: \spax(\CC)-\{0\} \to \sphere$ be the continuous
function given by $\proba(\var)=\var/\normcomp \var\normcomp $. Then
$\var= \normcomp \var\normcomp \proba(\var)$. Since we are on a
homogeneous affine variety, if $\var \in \spax (\CC)$, then also
$\compnumt \var \in \spax(\CC)$ for $\compnumt \in \CC$ (we need it
only for $\compnumt \in \RR_+$). A homogeneous function $\fuf$ on
$\spax (\CC)$ of degree $\degd$ has the property that $\fuf
(\compnumt \var)= \compnumt^\degd \fuf(\var)$ for $\compnumt \in
\CC^\times$. So in particular we have $\fuf(\var) = \fuf (\normcomp
\var\normcomp \frac{\var}{\normcomp \var\normcomp })= \normcomp
\var\normcomp ^\degd \fuf(\proba(\var))$ for $\var \neq 0$.

Define for $\var \in \spax (\CC)- \{0\}$ the functions $\contfu_i
(\var) = \normcomp \var\normcomp ^{\degd-\degd_\indi}
\contfuhi_\indi (\proba(\var))$. Since $\sphere$ is compact, the
continuous functions $ \contfuhi_\indi $ are bounded on $\sphere$
and since $\degd
> \degd_\indi$ we have $\lim_{\normcomp \var\normcomp  \to 0} \normcomp \var\normcomp ^{\degd-\degd_\indi} \contfuhi_\indi (\proba(\var))=0$
(in particular, this limit exists). So these functions extend to
continuous functions on $\spax (\CC)$. We get for $\var \neq 0$ the
equations
\begin{eqnarray*}
\fuf (\var) &=&\normcomp \var\normcomp ^\degd      \fuf( \proba
(\var) ) \cr &=& \normcomp \var\normcomp ^\degd \big( \contfuhi_1(
\proba(\var)) \fuf_1( \proba(\var)) \plusdots  \contfuhi_\numgen (
\proba(\var)) \fuf_\numgen  (\proba(\var)) \big)  \cr &=& \normcomp
\var\normcomp ^\degd    \big( \normcomp \var\normcomp ^{\degd_1-
\degd} \contfu_1(\var) \normcomp \var\normcomp ^{-\degd_1}
\fuf_1(\var) \plusdots \normcomp \var\normcomp ^{\degd_\numgen
-\degd} \contfu_\numgen (\var) \normcomp \var\normcomp
^{-\degd_\numgen } \fuf_\numgen (\var) \big) \cr &=& \contfu_1(\var)
\fuf_1 (\var)\plusdots \contfu_\numgen (\var)f_\numgen (\var) \, .
\end{eqnarray*}
By continuity, this equation holds also in the vertex.
\end{proof}

\ren{\runi}{{j}}

\begin{remark}
\label{contfuremark} We recall the argument that a family of
complex-valued functions $\fuf_1 \comdots \fuf_\numgen $ on a
topological space $\spax$ without common zeros generate the unit
ideal in $C_\CC(\spax)$. With $\fuf_\indfu $ also its complex
conjugate $\compconj{\fuf_\indfu }$ is continuous. Hence $\fug=
\sum_{\runi=1}^\numgen \compconj{\fuf_\runi } \fuf_\runi$ belongs to
the ideal $\ideal =(\fuf_1 \comdots \fuf_\numgen) C_\CC(\spax)$.
Since $\fug=\sum_{\runi=1}^\numgen \compconj{\fuf_\runi}f_\runi =
\sum_{\runi=1}^\numgen \normcomp \fuf_\runi\normcomp ^2 \geq 0$, the
set of zeros of this real valued function is the common zero locus
of the complex-valued functions $\fuf_1 \comdots \fuf_\numgen $. So
if the $\fuf_\indi $ have no common zero, then $\fug$ has no zero at
all and $1/\fug$ is a continuous function on $\spax$. Hence $\fug$
is a unit.

Therefore, in the situation of a primary ideal $\ideal =(\fuf_1
\comdots \fuf_\numgen  )$ in a standard-graded $\CC$-algebra $\ring$
and $\fuf \in \ring$, $\spax =\Spec \ring $, we get on $\spax
(\CC)-\{0\}$ the continuous coefficient functions
$$\contfuhi_\indfu = \frac{f \compconj{f}_\indfu}{f_1\compconj{f}_1 \plusdots \fuf_\numgen \compconj{f}_\numgen }
= \frac{f \compconj{f}_\indfu}{\normcomp \fuf_1\normcomp ^2
\plusdots \normcomp \fuf_\numgen\normcomp ^2 } \, .$$ If the
$\fuf_\indfu $ are homogeneous of the same degree, then the
real-homogeneous extension of the restriction to $\sphere$ used in
the proof of Theorem \ref{degreecriterioncont} gives us these
functions back. In general, if all these $\contfuhi_\indfu$ have a
limit in the origin, then they extend to continuous functions on
$\spax (\CC)$ and hence $\fuf \in (\fuf_1 \comdots
\fuf_\numgen)^\cont$. So these are natural candidates to look at for
continuous solutions. However, even if $\fuf \in (\fuf_1 \comdots
\fuf_\numgen)^\cont$ these functions do not always have a continuous
extension. For $\ideal =(\varx,\vary)$ and $\fuf = \varx$ we get
$\contfuhi_1= \frac{ \varx \compconj{\varx} }{\normcomp
\varx\normcomp ^2+\normcomp \vary\normcomp ^2}= \frac{\normcomp
\varx\normcomp ^2}{\normcomp \varx\normcomp ^2+\normcomp
\vary\normcomp ^2}$, which does not have a limit in the origin (look
at the limits for $\varx =0$ and $\vary =0$). The function is
however bounded around the origin, and in fact one can deduce from a
characterization of the integral closure due to Teissier that $\fuf$
belongs to the integral closure of $\fuf_1 \comdots \fuf_\numgen$ if
and only if one can write $\fuf= \fubound_1 \fuf_1 \plusdots
\fubound_\numgen \fuf_\numgen $ with locally bounded functions
$\fubound_\indfu$; see \cite{blicklebrennersubmersion}.
\end{remark}

\begin{corollary}
Let $\fuf \in \CC[\var_1 \comdots \var_\numvar]$ be a homogeneous
polynomial of degree $\degd \geq 1$. Then there exist continuous
functions $\contfu_\indfu:\CC^\numvar \to \CC$ such that $\fuf =
\contfu_1 \var_1^{\degd-1} \plusdots \contfu_\numvar
\var_\numvar^{\degd-1}$.
\end{corollary}
\begin{proof}
This follows directly from Theorem \ref{degreecriterioncont}.
\end{proof}

\ren{\dege}{{e}}

\begin{example}
Let $\ideal =(\var_1^{\degd_1} \comdots
\var_\numvar^{\degd_\numvar})$ and $\fuf  \in \CC[\var_1 \comdots
\var_\numvar]$ of degree $\degd > \max \degd_\indfu$. Then on
$\CC^\numvar-\{0\}$ we write $$\fuf (\var) = (\frac{\fuf(\var)
\compconj {\var}_1^{\degd_1}}{ \sum_{\runi=1}^\numvar \normcomp
\var_\runi\normcomp ^{2\degd_\runi} }) \var_1^{\degd_1} \plusdots
(\frac{\fuf (\var) \compconj{\var}_\numvar^{\degd_\numvar}}{
\sum_{\runi=1}^\numvar \normcomp \var_\runi\normcomp ^{2\degd_\runi}
}) \var_\numvar^{\degd_\numvar} \, .
$$
So in the proof of Theorem \ref{degreecriterioncont} we have
$\contfuhi_\indfu (\var) = \frac{\fuf(\var)
\compconj{\var}_\indfu^{\degd_\indfu}}{ \sum_{\runi =1}^\numvar
\normcomp \var_\runi\normcomp ^{2\degd_\runi} }$ and accordingly
$\contfu_\indfu (\var)= \normcomp \var\normcomp
^{\degd-\degd_\indfu} \contfuhi_\indfu( \frac{\var}{\normcomp
\var\normcomp }) = \normcomp \var\normcomp ^{\degd-\degd_\indfu}
\frac{ f(\frac{\var}{\normcomp \var\normcomp })
({\frac{\compconj{\var}_\indfu} {\normcomp \var\normcomp
}})^{\degd_\indfu} }{ \sum_{\runi =1}^\numvar \normcomp
\frac{\var_\runi}{\normcomp \var\normcomp }\normcomp ^{2
\degd_\runi}}$. If the $\degd_\indfu=\dege$ are constant, then this
is just $\contfu_\indfu (\var) = \frac{ \fuf (\var)
\compconj{\var}_\indfu^\dege} { \sum_{\runi =1}^\numvar \normcomp
\var_\runi\normcomp ^{2 \dege} }$. If $\ideal
=(\varx^\dege,\vary^\dege)$ and $\fuf =\varx^\expox \vary^\expoy$ is
a monomial, $\expox, \expoy < \dege$, then $\contfu_1(\varx, \vary)=
\frac{\varx^\expox \vary^\expoy \compconj{\varx}^\dege}{\normcomp
\varx\normcomp ^{2\dege}+\normcomp \vary\normcomp ^{2\dege} }
=\frac{\normcomp \varx\normcomp ^{2\dege} \vary^\expoy
\compconj{\varx}^{\dege - \expox}}{\normcomp \varx\normcomp
^{2\dege}+\normcomp \vary\normcomp ^{2\dege} } $ and
 $\contfu_2(\varx, \vary)
=\frac{\normcomp \vary\normcomp ^{2\dege} \varx^\expox
\compconj{\vary}^{\dege - \expoy}}{\normcomp \varx\normcomp
^{2\dege}+\normcomp \vary\normcomp ^{2\dege} } $.
\end{example}

\begin{remark}
\label{monomincontalt} We give an alternative construction for
continuous functions in the situation mentioned in the last example,
so $\ideal =(\varx^\expon,\vary^\expon) \subset \CC[\varx,\vary]$
and $\fuf = \varx^\expox\vary^\expoy$ with $\expox,\expoy <\expon$
and $\expox+\expoy >\expon$.

Consider the projective line $\PP^1_\CC \cong \sphere^2$ and denote
its complex variable by $\varz$, hence $\varz \in \CC$ or $\varz=
\infty$. Consider the function $\cofu (\varz) =
\frac{1}{\varz^\expon} - \frac{\varz^\expoy}{\varz^\expon}$ in a
(complex) neighborhood of $\varz =\infty$. It is a continuous
function in such a neighborhood (and even an algebraic function on
$\PP^1_\CC - \{0\}$). We can extend it as a continuous function to
the whole sphere, and denote the resulting function again by $ \cofu
(\varz)$. Note that the function $\varz^\expon \cofu (\varz) +
\varz^\expoy$ is also a continuous function on $\sphere^2$. This is
true in a neighborhood of $\infty$, because the function is even
constant $=1$ there, by the definition of $\cofu$ and for $\varz
\neq \infty$ it is clear anyway.

Now we use the substitution $\varz=\vary/\varx$. We set for
$(\varx,\vary) \neq (0,0)$
$$q_2 (\varx,\vary) =
\varx^{\expox+\expoy-\expon} \psi (\frac{\vary}{\varx}) \, .$$ As
$\fupsi$ is a bounded (since continuous) function on $\sphere^2$ and
$\expox+\expoy>\expon$ this function is continuous on $\CC^2-\{0\}$
and extends continuously to $\CC^2$ with value $0$ at the origin. We
set
$$\contfu_1(\varx,\vary) = - \varx^{\expox+\expoy-\expon}((\frac{\vary}{\varx})^\expon \psi ( \frac{\vary}{\varx})
+(\frac{\vary}{\varx})^\expoy) \, .$$ Again, this gives a continuous
function on $\CC^2$. We check for $\varx \neq 0$ (this we need for
resolving the bracket)
\begin{eqnarray*} \varx^\expon \contfu_1+\vary^\expon \contfu_2+\varx^\expox\vary^\expoy
&\!\!\!=\!\!\!& - \varx^\expon \varx^{\expox+\expoy-\expon}
\big(\!(\frac{\vary}{\varx})^\expon \fupsi ( \frac{\vary}{\varx})
+(\frac{\vary}{\varx})^\expoy\! \big) + \vary^\expon \varx^{\expox
+\expoy-\expon } \fupsi (\frac{\vary}{\varx})
+\varx^\expox\vary^\expoy \cr &\!\!\!=\!\!\!&
-\varx^{\expox+\expoy-\expon}\vary^\expon \fupsi (
\frac{\vary}{\varx})- \varx^{\expox}\vary^\expoy + \vary^\expon
\varx^{\expox+\expoy-\expon} \fupsi (\frac{\vary}{\varx})
+\varx^\expox \vary^\expoy \cr &\!\!\!=\!\!\!& 0 \, .
\end{eqnarray*}
For $\varx=0$ we have $\contfu_2 (0,\vary) = 0$, again since
$\expox+\expoy
>\expon$. Hence also under this condition the equation holds.
\end{remark}

\section{Powers of functions and powers of ideals}

We want to show that $\fuf^\expofu \in \ideal^\expoideal$ for
$\expofu < \expoideal $ implies that $\fuf \in \ideal^\cont$. For
this we first prove two lemmas.

\begin{lemma}
\label{persistence} The continuous closure is persistent, that is if
$\homring: \ring \to \rings$ is a $\CC$-algebra homomorphism of
finitely generated $\CC$-algebras, $\ideal$ is an ideal in $\ring$
and $\fuf \in \ring$ with $\fuf \in \ideal^\cont$, then also
$\homring(\fuf) \in (\ideal \rings)^\cont$.
\end{lemma}
\begin{proof}
Let $\ideal= (\runfuf)$. Then the condition means that $\fuf =
\contfu_1 \fuf_1 \plusdots \contfu_\numfu \fuf_\numfu$ with
continuous functions $\contfu_\indfu : (\Spec \ring)(\CC) \to \CC$.
Pulling these functions back via the continuous mapping $\homring^*:
(\Spec \rings)(\CC) \to (\Spec \ring)(\CC) $ gives $\homring (\fuf)
= ( \contfu_1 \circ \homring^*) \homring ( \fuf_1) \plusdots
(\contfu_\numfu \circ \homring^*) \homring (\fuf_\numfu)$. This
shows that $\homring(\fuf) \in (\homring ( \fuf_1) \comdots \homring
(\fuf_\numfu))^\cont = (\ideal \rings)^\cont$.
\end{proof}

The following Lemma deals with the generic situation of the proposed
question.

{ \ren{\multexp}{{\mu}}
\begin{lemma}
\label{potenzvariable} Let $ \expofu < \expoideal$ and let $$\ring =
\CC[ \varfu, \varcoef_\indcoef, |\indcoef|= \expoideal ,\vargen_1
\komdots \vargen_\numgen]/(\varfu^\expofu - \sum_{\indcoef,
|\indcoef|= \expoideal} \varcoef_\indcoef \vargen^\indcoef) \, .$$
Then $\varfu \in (\vargen_1 \komdots \vargen_\numgen)^\cont$ in
$\ring$.
\end{lemma}
\begin{proof}
According to Remark \ref{contfuremark} we look at the functions
$\frac{ \varfu \compconj{ \vargen}_\indgen}{ \normcomp
\vargen_1\normcomp ^2 \plusdots \normcomp \vargen_\numgen \normcomp
^2}$ defined on $D(\vargen_1 \komdots \vargen_\numgen) \subset
\spax(\CC)$ ($\spaxeqspecring$). We  want to show that they have
limit $0$ as a point converges to $V(\vargen_1 \komdots
\vargen_\numgen)$ and hence that they have a continuous extension
from $D(\vargen_1 \komdots \vargen_\numgen)$ to $\spax(\CC)$. The
modulus of such a function is $\frac{ \normcomp \varfu \normcomp
\normcomp \vargen_\indgen\normcomp }{ \normcomp  \vargen_1\normcomp
^2 \plusdots \normcomp \vargen_\numgen \normcomp ^2}
 = \frac{ \normcomp \varfu \normcomp  }{\sqrt{ \normcomp  \vargen_1\normcomp ^2 \plusdots
\normcomp \vargen_\numgen \normcomp ^2}} \frac{ \normcomp
\vargen_\indgen \normcomp }{\sqrt{ \normcomp  \vargen_1\normcomp ^2
\plusdots \normcomp \vargen_\numgen\normcomp ^2}} $ and the second
factor is bounded by $1$, so we only have to deal with the first
one. The function $\frac{ \normcomp \varfu \normcomp  }{\sqrt{
\normcomp \vargen_1\normcomp ^2 \plusdots \normcomp \vargen_\numgen
\normcomp ^2}} = \frac{ \normcomp \varfu \normcomp  }{ \normcomp
(\vargen_1 \komdots \vargen_\numgen )\normcomp }$ has limit $0$ for
a sequence $ \point_\indseq \in D(\vargen_1 \komdots
\vargen_\numgen)$ converging to $\point \in V(\vargen_1 \komdots
\vargen_\numgen)$ if and only if this is true for a (natural) power
of it. Hence we replace this function by $ \frac{ \normcomp \varfu
\normcomp ^\expofu }{ \normcomp ( \vargen_1 \komdots \vargen_\numgen
)\normcomp ^\expofu } = \frac{ \normcomp  \sum_\indcoef
\varcoef_\indcoef \vargen^\indcoef \normcomp  }{ \normcomp (
\vargen_1 \komdots \vargen_\numgen )\normcomp ^\expofu }$. We may
look at the summands separately. Since $\point_\indseq$ converges,
the values of $\varcoef_\indcoef (\point_\indseq)$ are bounded, and
so we consider $\frac{ \normcomp \vargen^\indcoef\normcomp
}{\normcomp (\vargen_1 \comdots \vargen_\numgen)\normcomp ^\expofu}
= \frac{\normcomp  \vargen_1^{\indcoef_1} \cdots
\vargen_\numgen^{\indcoef_\numgen}\normcomp }{\normcomp (\vargen_1
\comdots \vargen_\numgen)\normcomp ^\expofu}$, where
$|\indcoef|=\indcoef_1 \plusdots \indcoef_\numgen= \expoideal
> \expofu$. We can write $\frac{\normcomp  \vargen_1^{\indcoef_1}
\cdots \vargen_\numgen^{\indcoef_\numgen}\normcomp } {\normcomp
(\vargen_1 \comdots \vargen_\numgen)\normcomp ^\expofu} = $ \
$$\frac{\normcomp  \vargen_1^{\indcoef_1}\normcomp  \cdots
\normcomp \vargen_\numgen^{\indcoef_\numgen}\normcomp
}{\sqrt{\normcomp \vargen_1\normcomp ^2\! \plusdots \! \normcomp
\vargen_\numgen\normcomp ^2}^\expofu}= \vargen^\multexp
\big(\frac{\normcomp  \vargen_1 \normcomp }{\sqrt{\normcomp
\vargen_1\normcomp ^2 \! \plusdots\! \normcomp
\vargen_\numgen\normcomp ^2}}\big) ^{\expofu_1}
 \cdots \big(\frac{\normcomp  \vargen_\numgen \normcomp }{\sqrt{\normcomp \vargen_1\normcomp ^2
\! \plusdots\!  \normcomp \vargen_\numgen\normcomp ^2}}\big)
^{\expofu_\numgen} \, ,$$ where $|\multexp| = \expoideal-\expofu
>0$. As $\point_\indseq \to \point$, $\vargen_\indgen \to 0$ and so
this function goes to $0$.
\end{proof}}

\begin{theorem}
\label{powerscompare} Let $\ring$ be a finitely generated
$\CC$-algebra, let $\ideal$ be an ideal  and $\fuf \in \ring$.
Suppose that there exist numbers $\expofu < \expoideal $ such that
$\fuf^\expofu \in \ideal^\expoideal$. Then $\fuf \in \ideal^\cont$.
\end{theorem}
\begin{proof}
Let $\ideal=(\runfuf)$ and write the condition as $$\fuf^\expofu =
\sum_{\indcoef, |\indcoef|= \expoideal}   \fucoef_\indcoef
\fuf_1^{\indcoef_1}\cdots \fuf_\numgen^{\indcoef_\numgen}, \, \,
 \fucoef_\indcoef \in \ring \, .$$
We have a ring homomorphism from the generic situation to this
special situation, namely
$$  \CC[ \varfu, \varcoef_\indcoef, |\indcoef|= \expoideal ,\vargen_1
\komdots \vargen_\numgen]/(\varfu^\expofu - \sum_{\indcoef,
|\indcoef|= \expoideal} \varcoef_\indcoef \vargen^\indcoef) \longto
\ring, \,\,  \vargen_\indgen \mapsto \fuf_\indgen,\,
\varcoef_\indcoef \mapsto \fucoef _\indcoef,\, \varfu \mapsto \fuf
\,.
$$
Lemma \ref{potenzvariable} shows that $\varfu \in (\vargen_1
\comdots \vargen_\numgen)^\cont$ in the ring on the left hand side.
By the persistence of the continuous closure (Lemma
\ref{persistence})  $\fuf \in \ideal^\cont$ follows.
\end{proof}

\begin{remark}
\label{notweakenremark} The condition in Theorem \ref{powerscompare}
can not be weakened. Example \ref{weaknotaxesexample} below shows
that the condition $\fuf^\expofu \in \ideal^\expofu$ for all
$\expofu \geq 2$ does not imply $\fuf \in \ideal^\cont$. We will
show in Corollary \ref{orderbiggeragain} that under the assumption
of Theorem \ref{powerscompare} the element $\fuf$ belongs also to
the axes-closure.

Theorem \ref{powerscompare} has also a purely topological analog. If
$\spax$ is a topological space an $\fuf$ and $\runfuf$ are
continuous complex-valued functions on $\spax$ such that
$\fuf^\expofu \in (\runfuf)^\expoideal$ (in $C_\CC(\spax)$) for some
$\expofu < \expoideal$, then already $\fuf \in (\runfuf)$. We can
use an equation $\fuf^\expofu = \sum_{\indcoef, |\indcoef|=
\expoideal}   \fucoef_\indcoef \fuf_1^{\indcoef_1}\cdots
\fuf_\numgen^{\indcoef_\numgen}$ with $ \fucoef_\indcoef \in
C_\CC(\spax)$ to get a continuous mapping $(\fuf_1 \komdots
\fuf_\numgen,\fucoef_\indcoef ,\fuf): \spax \to \spav =
V(\varfu^\expofu - \sum_{\indcoef, |\indcoef|= \expoideal}
\varcoef_\indcoef \vargen^\indcoef) \subseteq \CC^\ell$. Since
$\varfu \in (\vargen_1 \comdots \vargen_\numgen)C_\CC(\spav)$ by
Lemma \ref{potenzvariable}, it follows $\fuf \in
(\runfuf)C_\CC(\spax)$.
\end{remark}

\section{$\field$-schemes of axes}
\label{axesschemesection}

We want to control the continuous closure algebraically, in
particular we want to find exclusion criteria for it. We fix a field
$\field$ and work in the category of finitely generated
$\field$-algebras. Whenever we talk about the continuous closure we
suppose that $\field=\CC$.

We have immediately the inclusions $\ideal \subseteq \ideal^\cont
\subseteq \rad (\ideal)$. In fact $\fuf  \in \rad (\ideal)$ holds if
and only if $\fuf = \fug_1\fuf_1 \plusdots \fug_\numgen
\fuf_\numgen$ with arbitrary functions $\fug_\indfu$. A better
algebraic approximation is given by the integral closure. One of the
equivalent characterizations of the integral closure
$\intclo{\ideal}$ of an ideal $\ideal $ in a $\field$-algebra
$\ring$ of finite type is that for every normal affine curve
$\schemetest=\Spec \ringt$ over $\field$ and every $\field$-morphism
$\homtest: \schemetest \to \Spec \ring $ we have $\varphi(\fuf) \in
\homtest(\ideal) \ringt$ \cite{blicklebrennersubmersion}. If
$\field=\CC$ and $\fuf  \in \ideal^\cont$, then for every such curve
it follows (by the persistence of the continuous closure, Lemma
\ref{persistence}) that $\varphi(\fuf) \in (\homtest
(\ideal)\ringt)^\cont$. But for a normal (=regular) curve always
$\ideal = \ideal^\cont$ holds, since $\ideal$ is locally a principal
ideal.

However, also the inclusion $\ideal ^\cont \subseteq
\intclo{\ideal}$ is strict. For example, $\varx \vary \in \intclo {
(\varx^2,\vary^2)}$, but $\varx \vary$ does not belong to the
continuous closure, as we will see (Corollary
\ref{monomialbordernotaxes}). Also, the weak subintegral closure
(\cite{leahyvitulliweaksubintegralclosure},
\cite{reidvitullimonomial}) is too big, as $\varx \vary^2 \in
(\varx^3,\varx^2 \vary, \vary^3)^\wsi$, but again not in the
continuous closure (Example \ref{weaknotaxesexample}). So our goal
here is to find a reasonably big category of one dimensional
$\field$ -schemes where the continuous closure is the identical
closure and which will be later on a test category for exclusion
results.

\begin{definition}
A \emph{scheme of axes} over $\field$ is a one-dimensional reduced
affine scheme $\schemeaxes = \Spec \ringt$ of finite type such that
its integral components are normal, meet transversally in one point
$\orig$ (the origin) and that the embedding dimension at $\orig$
equals the number of components. We call $\ringt$ a \emph{ring of
axes}.
\end{definition}

\ren{\ind}{i}

\begin{remark}
\label{axesremark} An equivalent characterization is that the
completion of the local ring at the origin $\orig$ is isomorphic to
$\fieldres(\orig)[[\varax_1 \comdots \varax_\numaxes ]]/(
\equataxes)$. A scheme of axes is a \emph{seminormal}
one-dimensional scheme and every seminormal one-dimensional local
scheme of finite type over an algebraically closed field looks after
completion like this (J. Tong called these singularities in a talk
in M\"unster, June 2006, `singularities of coordinate axe type'; for
the notion of seminormal rings see \cite{leahyvitulliseminormal},
\cite{leahyvitulliseminormalcorrection},
\cite{brennersemiintegraltest}). The easiest example is given by the
standard scheme of axes $\field[\varax_1 \comdots \varax_\numaxes
]/(\equataxes)$. The rings $\field[\varax, \varaxv]/(\varax \varaxv
(\varax + \varaxv)) $ and $\field[\varax, \varaxv]/(\varaxv (\varaxv
+ \varax^2)) $ do not give scheme of axes, though their components
are smooth and though they give for $\field=\CC$ topologically a
space of axes.

The property of being a scheme of axes depends on the field
$\field$. The ring $\ringt=\RR[\varx, \vary]/(\varx^2+\vary^2)$ is
seminormal, its normalization is $\CC[\varx]$ (with a non-trivial
extension of residue class fields). As $\Spec \ringt$ is integral
with one non-normal component, it is not a scheme of axes. The
tensoration $\ringt \tensor_\RR \CC=\CC[\varx, \vary]/(( \varx+i
\vary)(\varx- i \vary))$ is however a ring of axes. If $\ringt$ is a
$\field$-algebra such that $\ringt \tensor_\field \overline{\field}$
is a ring of axes, then there exist also a finite extension $\field
\subseteq \fieldl \subseteq  \overline{\field}$ such that $\ringt
\tensor_\field \fieldl$ is a ring of axes over $\fieldl$ and hence
over $\field$.

If $\schemeaxes_\indaxis$ are smooth affine curves over $\field$
with fixed $\field$-points $\orig_\indaxis \in
\schemeaxes_\indaxis$, $\indaxis = 1 \comdots \numaxes$, then
$$\schemeaxes =\{(\var_1 \comdots \var_\numaxes) \in \schemeaxes_1
\timesdots \schemeaxes_\numaxes : \exists \indaxis \mbox{ such that
} \var_\indaxisrun = \orig_\indaxisrun \mbox{ for all } \indaxisrun
\neq \indaxis \}
$$
is a scheme of axes with origin $\orig = \orig_1 \timesdots
\orig_\numaxes$ and components $\schemeaxes_\indaxis \cong \orig_1
\timesdots \orig_{\indaxis -1} \times  \schemeaxes_\indaxis \times
\orig_{\indaxis+1} \timesdots \orig_\numaxes $.
\end{remark}

Let $\schemeaxes=\Spec \ringt$ be a scheme of axes with axes
(integral components) $\axis_\indaxis  $, $\indaxis =1 \comdots
\numaxes$. After localizing at $\orig$, each $\axis_\indaxis  $ is
the spectrum of a discrete valuation domain. In particular, for
every $\indaxis $ there exists a valuation $\valaxis_\indaxis$ on
$\ringt$, which assigns to $\fuf \in \ringt$ the order
$\valaxis_\indaxis (\fuf)= \valaxis_\indaxis (\fuf
|{\axis_\indaxis})$ (which is $\infty$ if $\fuf
|{\axis_\indaxis}=0$). For an ideal $\ideal$ we set
$\ordaxis_\indaxis=\valaxis_\indaxis (\ideal) = \min
\{\valaxis_\indaxis (\fuf):\, \fuf \in \ideal \}$. We describe the
ideals in a complete scheme of axes.

\begin{lemma}
\label{axidealgenerators} Let $\field$  be a field and let
$\ringt=\field[[\varax_1 \komdots \varax_\numaxes ]]/(\equataxes)$
or $\ringt=\CC\brackconvl\varax_1 \komdots
\varax_\numaxes\brackconvr/(\equataxes)$, the ring of convergent
power series {\rm(}convergent in some open complex neighborhood of
the origin{\rm)}. Let $\ideal \subseteq \ringt$ be an ideal. Let
$\ordaxis _\indaxis $ be the order {\rm(}possibly $\infty${\rm)} of
the ideal on the $\indaxis$th axis given by $\varax_\indaxisrun=0$
for $\indaxisrun \neq \indaxis $. Then
$(\varax_\indaxis^{\ordaxis_\indaxis + 1} ) \subseteq \ideal $ for
the indices $\indaxis$ such that $\ordaxis_\indaxis $ is finite.
Moreover, the ideal $\ideal $ has after reindexing a system of
generators of the form
$$\fuf_\indgen = \varax_\indgen^{\ordaxis _\indgen} + \realnumlambda_{\indgen,\numgen + 1}
\varax_{\numgen +1}^{\ordaxis_{\numgen + 1}} \plusdots
\realnumlambda_{\indgen, \ello} \varax_{\ello}^{\ordaxis_{\ello}
},\, \, \indgen=1 \comdots \numgen \, .$$ Here
$\realnumlambda_{\indgen, \indaxis} \in \field$ and $\varax_\indaxis
$, $1 \leq \indaxis \leq \ello$, correspond to the components where
$\ordaxis _\indaxis $ is finite.
\end{lemma}
\begin{proof}
The statement is true for $\ideal =(1)$, so assume that $\ideal \neq
(1)$. Every non unit $\fuf \in \ringt$ can be written as $\fuf
=\sum_{\indaxis =1}^\numaxes \varax_\indaxis^{\degaxis_\indaxis}
\powser_\indaxis (\varax_\indaxis)$, where
$\degaxis_\indaxis=\valaxis_\indaxis(\fuf) \geq 1$ for
$\fuf|_{\axis_\indaxis} \neq 0$ and where $\powser_\indaxis =0$ or
$\powser_\indaxis$ is a (convergent) power series in
$\varax_\indaxis $ with non zero constant term (and a unit in
$\ringt$). If the order of $\ideal $ on the $\indaxis $th component
is $\ordaxis _\indaxis < \infty $, then there exists an element
$\fuf \in \ideal$ with $\varax_\indaxis$-exponent $\degaxis
_\indaxis =\ordaxis _\indaxis $. Then $\varax_\indaxis\fuf =
\varax_\indaxis^{\ordaxis_\indaxis+1} \powser_\indaxis
(\varax_\indaxis)$, and so $\varax_\indaxis^{\ordaxis_\indaxis+1}
\in \ideal$.

Assume after reindexing that the variables $\varax_1   \comdots
\varax_\ello$ occur with finite order in $\ideal $ (so the other
variables do not occur at all). We look at a system of generators
for $\ideal $ which encompasses the $\varax_\indaxis   ^{\ordaxis
_\indaxis +1}$, $\indaxis=1 \comdots \ello$. Then the other
generators can be written as $\fug = \sum_{\indaxis =1}^\ello
\realnumlambda_\indaxis (\fug) \varax_\indaxis ^{\ordaxis_\indaxis}$
with coefficients $\realnumlambda_\indaxis (\fug) \in \field$.
Applying Gauss elimination and reindexing gives a generating system
as stated, where also the $\varax_\indaxis ^{\ordaxis_\indaxis +1}$
are not necessary anymore.
\end{proof}

\begin{corollary}
\label{orderstrict} Let $\ringt=(\field[\varax_1   \comdots
\varax_\numaxes ]/(\equataxes))_{(\varax_1 \comdots \varax_\numaxes
)}$ or its completion $\ringt=\field[[\varax_1 \comdots
x_\numaxes]]/(\equataxes)$ and let $\ideal $ be an ideal, $\fuf \in
\ringt$. Let $\valaxis_\indaxis$ be the valuation on the
$\indaxis$th axes given by $\varax_\indaxisrun =0, \indaxisrun \neq
\indaxis $. If $\degaxis_\indaxis=\valaxis_\indaxis (\fuf)>
\ordaxis_\indaxis (\ideal)= \ordaxis_\indaxis$ for all $\indaxis $,
then $\fuf \in \ideal$.
\end{corollary}
\begin{proof}
We only deal with the complete case. By Lemma
\ref{axidealgenerators}, $\varax_\indaxis   ^{\ordaxis  _\indaxis
+1} \in \ideal$. We can write $\fuf = \varax_1 ^{\degaxis
_1}\powser_1 \plusdots \varax_\numaxes ^{\degaxis
_\numaxes}\powser_\numaxes$, where $\powser_\indaxis \in \ringt$.
The condition means that $\degaxis_\indaxis > \ordaxis _\indaxis $,
hence $\fuf \in (\varax_1^{\ordaxis_1+1} \comdots
\varax_\numaxes^{\ordaxis_\numaxes+1}) \subseteq \ideal$.
\end{proof}

We want to show now that for $\field = \CC$ the continuous closure
on a scheme of axes is the identical closure operation.

\begin{lemma}
\label{axescont} Let $\ringt= \CC\brackconvl\varax_1 \comdots
\varax_\numaxes \brackconvr/(\equataxes)$. Let $\ideal =(\fuf_1
\comdots \fuf_\numgen )$ and $\fuf \in \ringt$ be given. Suppose
that there exists a complex open neighborhood $\neigh$ of the origin
where $\fuf$ and $\fuf_\ind$ converge and that there exist
continuous complex-valued functions $\contfu_1 ,\ldots
,\contfu_\numgen $ defined on $\neigh$ such that $\fuf = \contfu_1
\fuf_1 \plusdots \contfu_\numgen \fuf_\numgen $ {\rm(}as continuous
functions on $\neigh${\rm)}. Then $\fuf  \in (\fuf_1 \comdots
\fuf_\numgen )$ holds in $\ringt$.
\end{lemma}
\begin{proof}
Let $\ideal  \neq (1)$ and let $\ordaxis  _\indaxis   $ be the order
of the ideal along the axis $\axis_\indaxis   $. As explained in
Lemma \ref{axidealgenerators}, there is a system of ideal generators
looking like
$$\fuf_\indfu = \varax_\indfu^{\ordaxis_\indfu} + \realnumlambda_{\indfu,\numgen
 +1} \varax_{\numgen
+1}^{\ordaxis_{\numgen+1}} \plusdots \realnumlambda_{\indfu,\ello}
\varax_{\ello}^{\ordaxis  _{\ello} },\, \, \indfu=1 \comdots \numgen
\, .$$ By the condition of continuity restricted to the axis
$\axis_\indaxis$ separately it follows that $\valaxis_\indaxis(\fuf)
\geq \ordaxis_\indaxis$. Hence we can write (after subtracting
elements in $\ideal$) $\fuf = \realnumbeta_{\numgen +1}
\varax_{\numgen +1}^{\ordaxis _{\numgen +1} } \plusdots
\realnumbeta_{\ello} \varax_{\ello}^{\ordaxis_{\ello}
}+\polp(\varax_{\ello+1}\comdots \varax_{\numaxes})$. Assume now
that the continuous functions $\contfu_\indj$ do their job, i.e.
$\sum_{\indj=1}^\numgen \contfu_\indj \fuf_\indj=\fuf$. It follows
immediately that $\polp(\varax_{\ello+1}\comdots \varax_{\numaxes})
=0$. Then on the axis $\axis_\indaxis   $, $ \indaxis    \leq
\numgen$ (given by $\varax_\indaxisrun = 0$ for $\indaxisrun \neq
\indaxis $), we get the condition
$$ \sum_{\indfu=1}^\numgen (\contfu_\indfu |_{\axis_\indaxis })
(\fuf_\indfu|_{\axis_\indaxis}) =(\contfu_\indaxis |_{\axis_\indaxis
}) \varax_\indaxis ^{\ordaxis_\indaxis   } =\fuf|_{\axis_\indaxis}=0
\, . $$ It follows that $ \contfu_\indaxis = 0$ on $\axis_\indaxis $
first outside the origin but by continuity also at the origin. So
$\contfu_\indfu(0)=0$ for all $\indfu=1 \comdots \numgen$.

Assume now that $\fuf \not\in \ideal$ and that
$\realnumbeta_\indaxis \neq 0$, $\indaxis \geq \numgen+1$. Then on
$\axis_\indaxis$ we get the equation $\sum_{\indfu=1}^\numgen
(\contfu_\indfu |_{\axis_\indaxis}) \realnumlambda_{\indfu,
\indaxis} \varax_\indaxis^{\ordaxis_\indaxis} =
\realnumbeta_\indaxis \varax_\indaxis   ^{ \ordaxis_\indaxis }$ and
so $\sum_{\indfu=1}^\numgen (\contfu_\indfu |_{\axis_\indaxis })
\realnumlambda_{\indfu, \indaxis}
 = \realnumbeta_\indaxis    $. This gives a contradiction at the origin. Hence
 $\realnumbeta_\indaxis    =0$ for $\indaxis > \numgen $ and so in fact $\fuf = 0$.
\end{proof}

\ren{\point}{{Q}}

\ren{\realnum}{r}

\begin{corollary}
\label{axescontident} Let $\schemeaxes= \Spec \ringt$ be a scheme of
axes over $\CC$, $\ideal \subseteq \ringt$ an ideal. Then
$\ideal^\cont =\ideal$.
\end{corollary}
\begin{proof}
Let $\fuf \in \ringt$ and $\fuf \in \ideal^\cont$. Then we have to
check that $\fuf \in \ideal \O_\point$ for every point $\point \in
\schemeaxes$. For $\point \neq \orig$ this is the (easier) case of
just one axis, so we may assume that $\point=\orig$ is the origin.
In a small neighborhood of $\orig$, $\schemeaxes$ is isomorphic as a
complex space to an open (ball) neighborhood of the origin of the
(standard) scheme of axes $\schemeaxesd$ in $\CC^\numaxes$. This is
also a homeomorphism, so we may assume that $\ideal=(\fuf_1 \comdots
\fuf_\numgen)$ and that $\fuf,\fuf_\indfu$ are holomorphic functions
defined on $\schemeaxesd_\realnum=\{\var \in \schemeaxesd :\,
\normcomp \var\normcomp < \realnum \}$. By assumption there exist
continuous functions $\contfu_\indfu$ such that
$\fuf=\contfu_1\fuf_1 \plusdots \contfu_\numgen \fuf_\numgen$ on
$\schemeaxesd_\realnum$. By Lemma \ref{axescont} this means that
$\fuf \in (\fuf_1 \comdots \fuf_\numgen)\CC\brackconvl\varax_1
\comdots \varax_\numaxes\brackconvr/(\equataxes)$. Hence also $\fuf
\in (\fuf_1 \comdots \fuf_\numgen)\CC[[\varax_1 \comdots
\varax_\numaxes]]/(\equataxes)$. This implies by faithful flatness
of the completion that $\fuf \in \ideal \ringt_\orig$.
\end{proof}

\section{The axes closure}
\label{axessection}

\begin{definition}
Let $\ideal$ be an ideal in a commutative ring $\ring$, $\fuf \in
\ring$. We say that $\fuf$ belongs to the ($\field$-)\emph{axes
closure}, $\fuf \in \ideal^\ax$ if for every ring homomorphism
$\homtest: \ring \to \ringt$, where $\ringt$ is a $\field$-scheme of
axes, we have $\homtest(\fuf) \in \homtest(\ideal)\ringt$.
\end{definition}

\begin{remark}
We will mainly be interested in the $(\field)$-axes closure for
ideals in the category of finitely generated $\field$-algebras. If
the base field is understood we will just talk about the axes
closure. Some of the results are independent of the base field, but
we will not treat this systematically.

The ($\field$)-axes closure of an ideal in a finitely generated
$\field$-algebra is inside the integral closure. The integral
closure can be tested by homomorphisms to discrete valuation
domains, but in the case of a $\field$-algebra $\ring$ of finite
type also by only looking at morphisms from affine normal
$\field$-curves to $\Spec \ring$. Every such curve is a scheme of
(one) axes. In \cite{brennersemiintegraltest} we show that the weak
subintegral closure as introduced by Leahy and Vitulli in
\cite{leahyvitulliweaksubintegralclosure} can be tested by morphisms
of scheme of axes with only two components (crosses). Hence also
$\ideal^\ax \subseteq \ideal^\wsi$.
\end{remark}

\begin{corollary}
\label{orderbigger} Let $\ring$ be a commutative ring, $\ideal
\subseteq \ring$ an ideal and $ \fuf \in \ring$. Suppose that for
every discrete valuation domain $(\dvd, \valu)$ and every ring
homomorphism $\homtest: \ring \to \dvd$ we have
$\valu(\homtest(\fuf))
> \valu(\homtest(\ideal))$. Then $\fuf \in \ideal^\ax$ {\rm(}independent of any $\field${\rm)}.
\end{corollary}
\begin{proof}
This follows directly from Corollary \ref{orderstrict} and the
definition of the axes closure.
\end{proof}

We can deduce the axes version of Theorem \ref{powerscompare}.

\begin{corollary}
\label{orderbiggeragain} Let $\ring$ be a commutative ring, $\ideal
\subseteq \ring$ an ideal and $ \fuf \in \ring$. Suppose that
$\fuf^\expofu \in \ideal^\expoideal$ for some $\expofu <
\expoideal$. Then $\fuf \in \ideal^\ax$.
\end{corollary}
\begin{proof}
We want to apply Corollary \ref{orderbigger} and look at $\homtest:
\ring \to \dvd$ to a discrete valuation domain. The containment $
\homtest(\fuf)^\expofu \in (\homtest(\ideal))^\expoideal$ implies
the relationship $\expofu \, \valu (\homtest(\fuf)) \geq
\expoideal\, \valu (\homtest(\ideal))$. This gives $\valu
(\homtest(\fuf))>\valu (\homtest(\ideal))$ for all valuations
$\valu$.
\end{proof}

We provide two further propositions which we will use in the
computation of the axes closure of a monomial ideal.

\ren{\schemeaxes}{{C}} \ren{\point}{{P}}

\begin{proposition}
\label{axesseminormal} Let $\field$ be an algebraically closed
field. Suppose that $\ring$ is a commutative ring and that $\fuf \in
\ideal^\ax$. Then $\homtest(\fuf) \in \homtest(\ideal) \ringt$ holds
also for every ring homomorphism to a seminormal one dimensional
affine ring over $\field$.
\end{proposition}
\begin{proof}
We may assume immediately that $\schemesn = \Spec \ring$ is an
affine seminormal one dimensional scheme of finite type over
$\field$, and we have to show that $\fuf \in \ideal$ under the
condition that this is true for every morphism $\schemeaxes \to
\schemesn$, where $\schemeaxes$ is a $\field$-scheme of axes. Since
the containment $\fuf \in \ideal$ is a local property we may assume
that $\schemesn$ has only one singular point. We know that the
completion at the singular point of $\schemesn$ looks like the
completion of a ring of axes. Let $\schemesn_1 \comdots
\schemesn_\numcomp$ be the integral components of $\schemesn $
(which are seminormal, but not normal in general), and let
$\schemeaxes_1 \comdots \schemeaxes_\numcomp$ be the normalizations.
On each $\schemeaxes_\ind$, let $\point_{\ind, \indj}$, $\indj \in
\Indj_\ind$, be the closed points mapping to the origin $\point_\ind
\in \schemesn_\ind$. Then the family of pointed curves
$\point_{\ind,\indj} \in \schemeaxes_\ind$, $\indj \in \Indj_\ind$,
$\ind=1 \comdots \numcomp$, glue together in the sense of Remark
\ref{axesremark} to get a scheme of axes over $\field$, say
$\schemeaxes$ and a morphism $\schemeaxes \to \schemesn$. This
morphism is in the completion at the singular point an isomorphism.
\end{proof}

\begin{lemma}
\label{finitepureaxes} Let $\field$ be an algebraically closed field
and let $\ringhom: \ring \to \ringsec$ be a finite and pure
{\rm(}e.g. faithfully flat{\rm)} homomorphism of $\field$-algebras.
Let $\ideal \subseteq \ring$ be an ideal and $\fuf \in \ring$.
Suppose that $\ringhom(\fuf) \in (\ringhom(\ideal) \ringsec)^\ax$.
Then already $\fuf \in \ideal^\ax$.
\end{lemma}
\begin{proof}
We may assume immediately that $\ring =\ringt$ is a ring of axes
over $\field$. Then $\ringsec$ is also a one dimensional
$\field$-algebra of finite type, let $\ringsec^\seminor$ be its
seminormalization. Since $\fuf \in (\ideal \ringsec)^\ax $ (we
denote the image of $\fuf$ in $\ringsec$ and in $\ringsec^\seminor$
again by $\fuf$), we know by Proposition \ref{axesseminormal} that
$\fuf \in \ideal \ringsec^\seminor$. The morphism $\Spec
\ringsec^\seminor \to \Spec \ringsec$ is a homeomorphism. By the
geometric criterion \cite{brennersemiintegraltest} for purety for
schemes over schemes of axes it follows that also $\ringt \to
\ringsec^\seminor$ is pure. Hence $\fuf \in \ideal $ in $\ringt$.
\end{proof}

In the case $\field=\CC$ we can easily establish the relation
between axes closure and continuous closure.

\begin{corollary}
\label{continaxes} Let $\ring$ be a $\CC$-algebra of finite type.
Then $\ideal ^\cont \subseteq \ideal^\ax$.
\end{corollary}
\begin{proof}
Let $\fuf \in \ideal^\cont$ and let $\homtest:\ring \to \ringt$ be a
homomorphism to a $\CC$-algebra of axes. By the persistence of the
continuous closure we have $\homtest(\fuf) \in (\homtest(\ideal)
\ringt)^\cont$. By Corollary \ref{axescontident} it follows that
$\homtest(\fuf) \in \homtest(\ideal) \ringt$, hence $\fuf \in
\ideal^\ax$.
\end{proof}

\begin{question}
\label{contaxquestion} Is for $\CC$-algebras of finite type the
continuous closure the same as the axes closure? Is this true for
the polynomial ring? A purely topological version of this question
is whether a continuous mapping $\spay \to \spax$ admits a
continuous section under the condition that for every continuous
mapping $\mortest:\spac \to \spax$ there exists a lifting
$\tilde{\mortest}:\spac \to \spay$, where $\spac$ is a space
consisting of a finite number of lines (or planes) meeting in one
point. Minimal topological requirements to make this a reasonable
question are that $\spax$ and $\spay$ are locally compact with only
finitely many components.
\end{question}

\section{The axes closure of a monomial ideal}
\label{axesmonomialsection}

We want to show that for a monomial ideal Question
\ref{contaxquestion} has a positive answer by computing explicitly
how its axes closure and continuous closure look like. Recall that
for a monomial ideal $\monomideal =(\var^\multexp_1 \comdots
\var^\multexp_\numgen)$ in $\field[\var_1 \comdots \var_\numvar]$ it
is helpful to consider the set of exponents $\setexp=\{ \multexp:
\var^\multexp \in \monomideal\}$ as a set of integral points
$\setexp \subseteq \NN^\numvar \subset \RR_+^\numvar$. For example,
the integral closure of a monomial ideal is given by the monomial
ideal consisting of $\var^\multexp$ such that $\multexp$ lies in the
convex hull $\intclo{\setexp}$ of $\setexp$ inside $\RR_+^\numvar$
\cite[Exercise 4.23]{eisenbud}. Also, the weak subintegral closure
of a monomial ideal has a combinatorial description in terms of the
exponents \cite[Proposition 3.3 and Theorem
4.11]{reidvitullimonomial}.

\begin{theorem}
\label{monomialaxes} Let $\ring=\field[\var_1 \comdots \var_\numvar
]$ be a polynomial ring and let $\ideal $ be a primary monomial
ideal. Then $\ideal ^\ax$ consists of $\ideal $ and of the monomials
in the interior of the convex hull $\intclo{\setexp}$, where
$\setexp=\{\multexp: \var^\multexp \in \ideal \}$ is the
corresponding set of exponents.
\end{theorem}
\begin{proof}
Suppose that a monomial is in the interior $\intclo{\setexp}^\circ $
of the complex hull. Then for every discrete valuation the monomial
has bigger order than the ideal (one can also show that the
situation of Corollary \ref{orderbiggeragain} is fulfilled), so by
Corollary \ref{orderbigger} it belongs to $\ideal^\ax$.

So let $\fuf$ be a polynomial consisting of monomials which lie on
the border and do not belong to the ideal itself. We may assume by
Remark \ref{axesremark} that $\field$ is algebraically closed. We
can enlarge the ideal by adding all monomials which are in the
support of this polynomial except one. Hence in particular we may
assume that $\fuf$ is a monomial, $\fuf= \var^\multexptau$,
$\multexptau \in \intclo{\setexp}- \intclo{\setexp}^\circ$,
$\multexptau \not \in \setexp$. Suppose that the exponent
$\multexptau$ lies on the affine hyperplane (forming a border of
$\setexp$) spanned by the exponents ${\multexpsigma_1} \comdots
{\multexpsigma_\numvar}$, $\multexpsigma_\indj \in \setexp$. A
suitable finite transformation $\var_\indj \mapsto
\var_\indj^{\degm_\indj}$ maps the monomials $\var^{\multexpsigma_1}
\comdots \var^{\multexpsigma_\numvar}$ to monomials of constant
degree $\mondeg$. All other monomials in $\setexp$ are mapped to
monomials of degree $\geq \mondeg$. By Lemma \ref{finitepureaxes},
the axes closure is not changed by such a finite free extension.

Hence we may assume by filling up the ideal that the ideal is given
by all monomials of fixed degree $\mondeg$ with the exception of one
monomial. Then Corollary \ref{monomialbordernotaxes} gives the
result.
\end{proof}



\ren{\indpoint}{{j}}

\synchronize{\indaxis}{\indpoint}

\ren{\indvar}{{r}}

\begin{proposition}
\label{propsamedegree} Let $\field$ be a field, $\ring=
\field[\var_1 \komdots \var_\numvar]$, $\degfu \in \NN$. Let $\ideal
\subset \ring$ be an ideal generated by polynomials $\runfu$ of
degree $\degfu$ and let $\fuf$ be another polynomial of degree
$\degfu$. Then $\fuf \in \ideal^\ax$ if and only if $\fuf \in
\ideal$.
\end{proposition}
\begin{proof}
We may immediately assume that $\field$  is an infinite field.
According to Lemma \ref{pointcheck}, let $\set$ be a set of
$\numaxes $ points in $\field^\numvar$ with the property that every
polynomial of degree $ \leq \degfu$ which vanishes on $\set$ must be
the zero polynomial. Let the $\numaxes$ points $\point_\indaxis$
have coordinates $\tupelpoint_\indaxis=(\tupelpoint_{\indaxis 1}
\comdots \tupelpoint_{\indaxis \numvar})$, $\indaxis =1 \comdots
\numaxes$. We consider the ring homomorphism
$$\homtest: \field[\var_1 \komdots \var_\numvar] \longrightarrow \field[\varax_1 \comdots \varax_\numaxes]/(\equataxes)$$
given by
$$\var_\indvar \longmapsto \tupelpoint_{1\indvar} \varax_1 \plusdots \tupelpoint_{\numaxes \indvar} \varax_\numaxes \, .$$
A polynomial $\fug$ of degree $\degfu$ is sent to $\homtest(\fug)=
\fug(\tupelpoint_1) \varax_1^\degfu \plusdots
\fug(\tupelpoint_\numaxes) \varax_\numaxes^\degfu$, where
$\fug(\tupelpoint_\indpoint) = \fug( \tupelpoint_{\indpoint 1}
\comdots \tupelpoint_{\indpoint \numvar})$. Hence $\homtest (\fuf)
\in (\homtest(\fuf_1) \comdots \homtest(\fuf_\numfu))$ if and only
if there exist $ \coefc_\indfu \in \field$, $\indfu =1 \comdots
\numfu$, such that $\homtest(\fuf) = \sum_{\indfu=1}^\numfu
\coefc_\indfu \homtest(\fuf_\indfu)$, and this means by looking at
each component (each coefficient of $\varax_\indaxis^\degfu$) that
$\fuf(\tupelpoint_\indpoint) = \sum_{\indfu=1}^\numfu \coefc_\indfu
\fuf_\indfu(\tupelpoint_\indpoint)$ for all points
$\tupelpoint_\indpoint$, $\indpoint = 1 \comdots \numaxes$. This
means by the choice of $\tupelpoint_\indpoint$ that $\fuf =
\sum_{\indfu=1}^\numfu \coefc_\indfu \fuf_\indfu$.
\end{proof}

\begin{corollary}
\label{monomialbordernotaxes} Let $\field$ be a field, $\ring=
\field[\var_1 \komdots \var_\numvar]$, $\mondeg \in \NN$. Let
$\ideal  \subset \ring$ be a monomial ideal generated by all
monomials of degree $\mondeg$ with the exception of
$\var^\expot=\var_1^{\expot_1} \cdots \var_\numvar
^{\expot_\numvar}$, $\sum_{\indj =1}^\numvar \expot_\indj =
\mondeg$. Then $\var^\expot \not \in \ideal^\ax$.
\end{corollary}
\begin{proof}
This is a special case of Proposition \ref{propsamedegree}
\end{proof}

\ren{\dege}{{d}}

\ren{\dime}{{m}}

\begin{lemma}
\label{pointcheck} Let $\field$ be an infinite field, let $\dege,
\dime \in \NN$. Then there exist $\numpoints$ points in
$\field^\dime$ such that every polynomial of degree $\leq \dege$
which vanishes at all these points is the zero polynomial.
\end{lemma}
\begin{proof}
This can be proved by induction over $\dime$, see \cite[Satz
54.7]{schejastorch2}.
\end{proof}

\begin{example}
\label{weaknotaxesexample} Even for monomial ideals in a
two-dimensional polynomial ring the axes closure (and the continuous
closure) is smaller than the weak subintegral closure. The easiest
example is the ideal $\ideal =(\varx^3,\varx^2\vary,\vary^3)
\subseteq \ring=\field[\varx, \vary]$ and $\fuf = \varx \vary^2$.
Then by \cite[Proposition 3.3 and Example 4.12]{reidvitullimonomial}
we have $\varx \vary^2 \in \ideal^\wsi$. We have directly $(\varx
\vary^2)^\expofu \in \ideal^\expofu$ for $\expofu \geq 2$. Consider
according to Proposition \ref{propsamedegree} the homomorphism
$\ring \to \ringt= \field[\varax_1,\varax_2,\varax_3,
\varax_4]/(\equataxes)$ given by ($\Char (\field) \neq 2$)
$$\varx \mapsto \fug= \varax_1+\varax_2 +\varax_4 \mbox{  and } \vary \mapsto \fuh=\varax_1+\varax_3+2 \varax_4
\, .$$ Then in $\ringt$ we have
$$\fug^3 = \varax_1^3+\varax_2^3 +\varax_4^3,\, \, \, \fug^2\fuh = \varax_1^3 +2\varax_4^3, \, \, \,\fuh^3=\varax
_1^3+\varax_3^3+8\varax_4^3 $$ and $\fug\fuh^2 =\varax_1^3
+4\varax_4^3$. Looking at the coefficients we see that $\fug\fuh^2$
is not in the extended ideal $(\fug^3,\fug^2\fuh,\fuh^3)$. Hence
$\varx\vary^2 \not\in \ideal^\ax$. It follows also for $\field=\CC$
that $\varx \vary^2 \not\in \ideal^\cont$.
\end{example}

\ren{\realnum}{\alpha}

\section{The continuous closure of a monomial ideal}
\label{contmonomialsection}

We can now put the previous results together and compute the
continuous closure of a monomial ideal.

\begin{theorem}
\label{monomialaxescont} Let $\ring=\CC[\var_1  \komdots
\var_\numvar ]$ and $\monomideal  \subset \ring$ be a primary
monomial ideal given by the monomials $\var^\multexp$, $\multexp \in
\setexp$. Then $\monomideal^\cont = \monomideal^\ax$, and this is
the monomial ideal given by the monomials in $\ideal$ and the
monomials with exponents in the interior of the convex hull
$\intclo{\setexp}$ of $\setexp$ in $\RR_+^\numvar$.
\end{theorem}
\begin{proof}
We have shown in Corollary \ref{continaxes} that $\monomideal ^\cont
\subseteq \monomideal^\ax$ and in Theorem \ref{monomialaxes} that
$\monomideal ^\ax$ has this description. So we only have to show
that the monomials with exponent in the interior
$\intclo{\setexp}^\circ$ are inside the continuous closure. If
$\var^\multexptau$ is inside the interior, then there exist
monomials $\var^{\multexpsigma_1} \comdots
\var^{\multexpsigma_\numvar} \in \monomideal$ such that the exponent
$\multexptau $ lies ``above'' the affine hyperplane spanned by the
exponents $ \multexpsigma_1 \comdots \multexpsigma_\numvar$ of these
monomials. Let $\multexpsigma_\indi= \multexpsigma_{\indi,\indj}$.
There exist positive numbers $\coeflins_1 \comdots \coeflins_\numvar
\in \NN$ such that $\sum_{\indj=1}^\numvar \coeflins_\indj
\realnum_\indj = \mondeg \in \NN $ if and only if $\realnum =
(\realnum_1 \comdots \realnum_\numvar )$ belongs to the affine
hyperplane spanned by the monomials. In particular,
$\sum_{\indj=1}^\numvar \coeflins_\indj \multexpsigma_{\indi,\indj}=
\mondeg$ for $\indi=1 \comdots \numvar$ and $\sum_{\indj=1}^\numvar
\coeflins_\indj \multexptau_{\indi,\indj} > \mondeg $. Applying the
ring homomorphism $\var_\indj \mapsto \var_\indj^{\coeflins_\indj}$
brings our situation to an equivalent situation (by Lemma
\ref{finiteextension}) where the generating monomials have the same
degree $\mondeg$ and the monomial in question has strictly bigger
degree. Hence by Theorem \ref{degreecriterioncont} it belongs to the
continuous closure of the ideal (or one can show after this step
that $(\varx^\multexptau)^\expofu \in \ideal^\expoideal$, $\expofu <
\expoideal$, and then apply Theorem \ref{powerscompare}).
\end{proof}

\ren{\ringr}{{S}}

\begin{lemma}
\label{finiteextension} Let $ \contmap :\CC^\numvar \to \CC^\numvar$
be
given by $(\var_1 \komdots \var_\numvar) \to
(\var_1^{\degm_1} \komdots \var_\numvar ^{\degm_\numvar})$. Let
$\ideal =(\fuf_1 \comdots \fuf_\numgen)$ be an ideal in the ring
$\ringr =C_\CC (\CC^\numvar)$ of continuous complex-valued functions
on $\CC^\numvar$. Then for $\fuf \in \ringr$ we have $\fuf \in
\ideal$ if and only if $\fuf \circ \contmap  \in (\ideal \circ
\contmap )\ringr $.
\end{lemma}
\begin{proof}
The ring extension $\ringr \subseteq \ringr$ given by $\fuf \mapsto
\fuf \circ \contmap$ is a direct summand exactly as in the case of a
polynomial ring because of the existence of the trace map. The
Galois group $\group=\ZZ/\degm_1 \timesdots \ZZ/\degm_\numvar$ acts
on $\CC^\numvar$ (and on $\ringr$) with quotient $\CC^\numvar$ (and
invariant ring $\ringr$). An element $\contfu \in \ringr$ yields the
invariant element $\sum_{\elgroup \in \group} \contfu \elgroup$,
where $(\contfu \elgroup)(\var) =\contfu( \elgroup_1\var_1 \comdots
\elgroup_\numvar \var_\numvar)$. Applying this to an equation $ \fuf
\circ \contmap = \contfu_1 (\fuf_1 \circ \contmap ) \plusdots
\contfu_\numgen (\fuf_\numgen \circ \contmap )$ gives a similar
equation for $\fuf$ and $\fuf_1 \comdots \fuf_\numgen $.
\end{proof}

\begin{example}
The ideal $\ideal =(\varx^2,\vary^5)$ in $\CC[\varx,\vary]$ and
$\fuf = \varx\vary^3$ gives an example where $\fuf$ does not fulfill
the degree assumption of Theorem \ref{degreecriterioncont}, but it
is due to Theorem \ref{monomialaxescont} in the continuous closure
anyway.
\end{example}

\begin{question}
Is it possible to give an effective criterion for $\fuf \in (\fuf_1
\comdots \fuf_\numgen)^\cont$ (in a polynomial ring)? Is, if the
ideal generators form a Gr\"obner basis, the containment inside the
continuous closure just a question of whether the initial term of
$\fuf$ belongs (in an inductive sense) to the continuous closure of
the initial ideal.
\end{question}

\begin{question}
Is the tight closure of an ideal always inside the continuous
closure? This is for normal $\CC$-domains of finite type a
reasonable question. The continuous closure is often bigger, as
already the regular case shows. In the homogeneous parameter case
$\ideal =(\fuf_1 \comdots \fuf_\numgen )$ in dimension $\geq 2$, we
have by the Theorem of Hara \cite[Theorem 6.1]{hunekeparameter} that
$\fuf \in \ideal^*$ implies (if $\fuf \not\in \ideal$) $\deg(\fuf)
\geq \sum_{\runi =1}^\numgen \deg(\fuf_\runi)$ and so in particular
$\deg(\fuf)
> \deg(\fuf_\runi)$ for all $\runi$, hence $\fuf \in \ideal^\cont$.

There is no inclusion between the continuous closure and the regular
closure. It is known that $\varzz \in (\varxx, \varyy)^{\reg}$ in
$\ring=\field[\varxx,\varyy,\varzz]/(\varzz^3 - \varxx^3 -
\varyy^3)$. However, $\ring/(\varxx, (\varyy - \varzz)(\varyy-
\unitroot_3 \varzz)) \cong \field[\varyy, \varzz]/ (\varyy -
\varzz)(\varyy- \unitroot_3 \varzz)$ is a ring of two axes, but
$\varzz $ is not a multiple of $\varyy$ in this ring.
\end{question}

\begin{remark}
The results of this paper hold probably also if we replace
$C_\CC(\spax)$ by $C_\CC^\infty(\spax)$, the ring of smooth
functions i.e. functions which are differentiable in the real sense
of arbitrary order. However,  modifications are needed as the
continuous solutions in Theorem \ref{degreecriterioncont} and
Theorem \ref{powerscompare} are not smooth. The functions used in
Remark \ref{monomincontalt} can be easily replaced by smooth
functions.
\end{remark}

\bibliographystyle{plain}

\end{document}